# Implementing the Law of Sines to solve SAS triangles

Konstantine Zelator
Dept. of Math and Computer Science[1]
Rhode Island College
600 Mount Pleasant Avenue
Providence, RI 02908
U.S.A.
e-mail : kzelator@ric.edu
konstantine_zelator@yahoo.com

## 1. Introduction

Most trigonometry or precalculus texts offer a fairy brief treatment of the subject of solving a triangle; while only a few contain a more extensive analysis of the subject matter. By "solving a triangle" , we refer to the problem of determining all three sidelengths and interior triangle angle on given information on some of the triangle's angles and/or sidelengths. Depending on the specific information, a unique triangle may be formed, more than one triangle, or no triangle at all. Robert Blitzer's book, *Precalculus,* is among those books that contain a typical treatment of this subject (see [1]). Now, the case *SAS* refers to the situation wherein two sidelengths are given, as well as the degree measure of the angle contained between the two sides. As we know from Euclidean geometry, such a triangle is uniquely determined; meaning that any two triangles constructed with these specifications must be congruent. It appears that not only Blitzer's book, but universally, in any book or text with the *SAS* case, the Law of Cosines is employed. For example, let us say that the lengths $a$ and $b$ are given; as well as the degree measure $\omega$ of the angle $\angle BCA$ (see Figure 1). Then, by using the Law of Cosines, one calculates the value of $c$ , the length of $\overline{BA}$ ; and by using the Law of Sines, one determines the values of $\sin\varphi$ and $\sin\theta$. By the use of the inverse function on a calculator if necessary, one determines the degree measures $\varphi$ and $\theta$.

[1]Effective August 1, 2009 and for the academic year 2009-2010:
Konstantine Zelator
Department of Mathematics
301 Thackeray Hall
139 University Place
University of Pittsburgh
Pittsburgh, PA 15260
e-mail: kzet159@pitt.edu

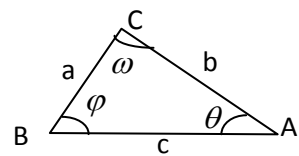

**Figure 1**

In this work, we present an alternative approach to solving an *SAS* triangle. An approach whose key ingredient is the Law of Sines.

*Law of Sines* 
$$\left\{\frac{a}{\sin\theta} = \frac{b}{\sin\varphi} = \frac{c}{\sin\omega}\right\} \qquad (1)$$

## 2      Basic trigonometric identities

The following identities are widely very well known and can be found in any standard trigonometry or precalculus text.

For any angle degree measures $\alpha$ and $\beta$,

$$\left\{\begin{array}{l}\cos(\alpha+\beta) = \cos\alpha\cos\beta - \sin\alpha\sin\beta \\ \sin(\alpha+\beta) = \sin\alpha\cos\beta + \cos\alpha\sin\beta\end{array}\right\} \qquad (2)$$

For any degree measures $\alpha$ and $\beta$,

$$\left\{\begin{array}{l}\sin\alpha + \sin\beta = 2\sin\left(\frac{\alpha+\beta}{2}\right)\cos\left(\frac{\alpha-\beta}{2}\right) \\ \sin\alpha - \sin\beta = 2\sin\left(\frac{\alpha-\beta}{2}\right)\cos\left(\frac{\alpha+\beta}{2}\right)\end{array}\right\} \qquad (3)$$

For any angle degree measure $\alpha$ **not** of the form $180k$ or $90+180k$; $\alpha \neq 180k, 90+180k$ ($k$ *any* integer),

$$\left\{\begin{array}{l}\tan(90-\alpha) = \cot\alpha \\ \cot(90-\alpha) = \tan\alpha\end{array}\right\}. \qquad (4)$$

For any angle degree measure $\alpha \neq 360k+180$ and $\alpha \neq 180k+90$ ($k$ *any* integer),

$$\tan \alpha = \frac{2 \tan\left(\frac{\alpha}{2}\right)}{1 - \tan^2\left(\frac{\alpha}{2}\right)}. \tag{5}$$

## 3      A Lemma and its proof

**Lemma 3.1.** Suppose that $a, b, c, d$ are real numbers such that $bd \neq 0$, $a + b \neq 0$, and $c + d \neq 0$. ( *In particular, this hypothesis is satisfied when $a, b, c, d$ are positive.*)

Then, $\frac{a}{b} = \frac{c}{d}$ if, and only if, $\frac{a-b}{a+b} = \frac{c-d}{c+d}$.

Proof. Suppose that $\frac{a}{b} = \frac{c}{d} = k$. Then $a = bk$ and $c = dk$. We have

$$\frac{a-b}{a+b} = \frac{bk-b}{bk+b} = \frac{b(k-1)}{b(k+1)} = \frac{k-1}{k+1} \tag{6}$$

and also

$$\frac{c-d}{c+d} = \frac{d\left(\frac{c}{d}-1\right)}{d\left(\frac{c}{d}+1\right)} = \frac{k-1}{k+1}. \tag{7}$$

From (6) and (7), it follows that $\frac{a-b}{a+b} = \frac{c-d}{c+d}$.

Conversely, assume that $\frac{a-b}{a+b} = \frac{c-d}{c+d}$, which implies

$(a-b)(c+d) = (a+b)(c-d) \Leftrightarrow ac + ad - bc - bd = ac - ad + bc - bd$, and thus, $2bc = 2ad$; and since $bd \neq 0$, we obtain $\frac{2bc}{bd} = \frac{2ad}{bd}$; $\frac{a}{b} = \frac{c}{d}$.

## 4      Solving an SAS triangle using the Law of Sines

Suppose that the lengths $a$ and $b$, as well as the degree measure $\omega$, are given in Figure 1. From the Law of Sines (1) we obtain

$$\frac{a}{b} = \frac{\sin\theta}{\sin\varphi} \tag{8}$$

Since $\theta$ and $\varphi$ are degree measures of triangle angles, we have $0 < \theta, \varphi < 180$, and so $\sin\theta > 0$ and $\sin\varphi > 0$; and, of course, $a > 0, b > 0$. By Lemma 1, (8) implies

$\dfrac{\sin\theta - \sin\varphi}{\sin\theta + \sin\varphi} = \dfrac{a-b}{a+b}$, and by identities (3),

$$\frac{2\sin\left(\frac{\theta-\varphi}{2}\right)\cos\left(\frac{\theta+\varphi}{2}\right)}{2\sin\left(\frac{\theta+\varphi}{2}\right)\cos\left(\frac{\theta-\varphi}{2}\right)} = \frac{a-b}{a+b} \tag{9}$$

Since $\theta, \varphi, \omega$ are triangle angle degree measures, the following hold true: $\theta + \varphi + \omega = 180$ and

$0 < \theta, \varphi, \omega < 180$. From which we obtain $\theta + \varphi = 180 - \omega$; $\dfrac{\theta+\varphi}{2} = 90 - \dfrac{\omega}{2}$ and

$0 < \omega < 180 \Leftrightarrow 0 < \dfrac{\omega}{2} < 90 \Leftrightarrow -90 < -\dfrac{\omega}{2} < 0 \Leftrightarrow 0 < 90 - \dfrac{\omega}{2} < 90$. This then shows that the conditions of identity (5) are satisfied and therefore,

$$\tan\left(\frac{\omega}{2}\right) = \cot\left(\frac{\theta+\varphi}{2}\right) \text{ and } \cot\left(\frac{\omega}{2}\right) = \tan\left(\frac{\theta+\varphi}{2}\right) \tag{10}$$

Also, note that form $-90 < -\dfrac{\varphi}{2} < 0$ and $0 < \dfrac{\theta}{2} < 90$ one easily obtains

$$-90 < \frac{\theta-\varphi}{2} < 90. \tag{11}$$

From (9) and (10) we infer that,

$$\boxed{\tan\left(\frac{\theta-\varphi}{2}\right) = \cot\left(\frac{\omega}{2}\right) \cdot \left(\frac{a-b}{a+b}\right)} \tag{12}$$

A well-known fact from trigonometry is that for every real number $v$, there is a unique degree measure $t$ such that $-90 < t < 90$ and $\tan t = v$. In other words, $t = \arctan v$ or $t = \tan^{-1} v$ with $t$ measured in degrees. Let us rewrite

$$-90 < t < 90 \text{ and } \tan t = v \tag{13}$$

Now, given $a, b,$ and $\omega$, we first calculate the value $v$ of the right-hand side of (12):

$$v = \cot\left(\frac{\omega}{2}\right) \cdot \left(\frac{a-b}{a+b}\right) \qquad (14)$$

$$v = \frac{1}{\tan\left(\frac{\omega}{2}\right)} \cdot \left(\frac{a-b}{a+b}\right)$$

Note that a calculator only has the $\boxed{\tan}$ key, which means that we cannot calculate the value of $\cot\left(\frac{\omega}{2}\right)$ directly.

Now, let $t$ be the unique degree measure such that

$\tan t = v$ and with $-90 < t < 90$.

Then, by (11), (12) , (13), and (14) it follows that

$$\left.\begin{array}{l} t = \dfrac{\theta - \varphi}{2}; \theta - \varphi = 2t. \\ \text{and also } \theta + \varphi = 180 - \omega \end{array}\right\} \qquad (15)$$

From (15) we obtain the formulas,

$$\boxed{\theta = 90 + t - \frac{\omega}{2}, \varphi = 90 - \left(t + \frac{\omega}{2}\right)} \qquad (16)$$

Remember that after we compute the value of $v$ in (14) by putting our calculator in degree mode, we then use the $\boxed{\tan^{-1}}$ key to find the value of $t$; and by the use of (16), we determine the degree values of $\theta$ and $\varphi$. Finally, from the Law of Sines (1) we can find the value of the sidelength $c$ by using either of $c = \dfrac{a \sin \omega}{\sin \theta}$ or $c = \dfrac{b \sin \omega}{\sin \phi}$.

## 5 Examples

1. Solve the *SAS* triangle that satisfies $\tan\left(\dfrac{\omega}{2}\right) = 2, a = 3,$ and $b = 1$.

*Solution:* In this example, the degree measure $\omega$ is not given directly. By (5) , we have $\tan \omega = \dfrac{4}{1-4} = -\dfrac{4}{3}$. This tells us that since $\omega$ is a triangle angle, $90 < \omega < 180$. By pressing that $\boxed{\tan^{-1}}$ key on key on the calculator in the degree mode, we obtain $\tan^{-1}\left(-\dfrac{4}{3}\right) \approx -53.13010235$. Thus, to obtain the (approximate) value of $\omega$, we must add $180$ to the answer obtained: $\boxed{\omega \approx 126.8698976}$.

Next we compute the value of $v$ in (14) with $\cot\left(\dfrac{\omega}{2}\right) = \dfrac{1}{2}, a = 3,$ and $b = 1.$ We find $v = \dfrac{1}{4}$ and so,

$$t = \tan^{-1}\left(\dfrac{1}{4}\right) \approx 14.03624347.$$

Applying formulas (16) gives $\boxed{\theta \approx 40.60129465 \text{ and } \varphi = 12.52880771}$.

Two decimal place approximations for the above angles are $\omega \approx 126.87, \theta \approx 40.60, \varphi \approx 12.53$ (with the sum being exactly equal to 180). Using the calculator again, we find that $\sin \omega \approx 0.8$ and $\sin \theta \approx 0.650791373.$ Hence,

$$c = \dfrac{a \sin \omega}{\sin \theta} \approx \dfrac{(3)(0.8)}{0.650791373} \approx \boxed{3.687817786}$$

2. What precise condition must the angle $\omega$, and the lengths $a$ and $b$ satisfy if $t = 45$?

   *Solution*: If $t = 45$, then $v = \tan t = \tan 45 = 1,$ and so by (14) we obtain $\tan\left(\dfrac{\omega}{2}\right) = \dfrac{a-b}{a+b}.$

   If we substitute for $\tan\left(\dfrac{\omega}{2}\right)$ in (5) we obtain, after some basic algebra $\boxed{\tan \omega = \dfrac{a^2 - b^2}{2ab}}$,

   is the precise condition sought after.

3. Suppose that $r$ is a given number greater than 1; and that $\dfrac{a}{b} = r$. In other words, we require that the sidelengths $a$ and $b$ have a prescribed ratio $r > 1$. Also, let $\omega = 60.$

   (i) Find the exact values of $\sin \theta$ and $\sin \varphi$ in terms of $r$.

   (ii) Express the third sidelength $c$ in terms of $r$ and $b$.

   (iii) Find the exact values of $\theta$ and $\varphi$ when $r = 2$.

   (iv) What can be said about the angles $\theta$ and $\varphi$, if we allow $r$ to vary in the open interval $(1, +\infty)$?

*Solution:*

(i) From $a = br$ and $\cot\left(\dfrac{\omega}{2}\right) = \cot(30) = \sqrt{3}$ we obtain $v = \dfrac{\sqrt{3} \cdot b \cdot (r-1)}{b(r+1)}$;

Any by (13) we arrive at $\tan t = \sqrt{3}\left(\dfrac{r-1}{r+1}\right) > 0$, since $r > 1$. Since $-90 < t < 90$ and $\tan t > 0$, it follows that $0 < t < 90$, which implies $\sin t > 0$ and $\cos t > 0$. We apply the identity $\sec^2 t = \tan^2 t + 1$; which implies, since $\cos t > 0$ and $\sec t = \dfrac{1}{\cos t}$, that

$$\cos t = \dfrac{1}{\sqrt{\tan^2 t + 1}} = \dfrac{1}{\sqrt{3\left(\dfrac{r-1}{r+1}\right)^2 + 1}} \underset{\text{algebra}}{=} \dfrac{r+1}{\sqrt{3(r-1)^2 + (r+1)^2}}$$

Next, $\sin t = (\cos t)(\tan t) = \left[\dfrac{r+1}{\sqrt{3(r-1)^2 + (r+1)^2}}\right] \cdot \sqrt{3}\left(\dfrac{r-1}{r+1}\right) = \dfrac{\sqrt{3}(r-1)}{\sqrt{3(r-1)^2 + (r+1)^2}}$.

Now apply the well known identity $\sin(90 + \alpha) = \cos \alpha$ with $\alpha = t - \dfrac{\omega}{2}$.

By (16), we have $\sin\theta = \cos\left(t - \dfrac{\omega}{2}\right) = \cos t \cos\left(\dfrac{\omega}{2}\right) + \sin t \sin\left(\dfrac{\omega}{2}\right)$; and by putting

$\cos\left(\dfrac{\omega}{2}\right) = \cos(30) = \dfrac{\sqrt{3}}{2}, \sin\left(\dfrac{\omega}{2}\right) = \sin(30) = \dfrac{1}{2}$, and employing the expressions we obtained above for $\sin t$ and $\cos t$ in terms of $r$, we get (after some careful algebra)

$$\boxed{\sin\theta = \dfrac{r\sqrt{3}}{\sqrt{3(r-1)^2 + (r+1)^2}}} = \dfrac{r\sqrt{3}}{2\sqrt{r^2 - r + 1}}.$$

From the Law of Sines (1) we also have $\dfrac{\sin\theta}{\sin\varphi} = \dfrac{a}{b} = r$; or equivalently, $\sin\varphi = \dfrac{1}{r} \cdot \sin\theta \Rightarrow$ (by the above formula)

$$\boxed{\sin\phi = \dfrac{\sqrt{3}}{\sqrt{3(r-1)^2 + (r+1)^2}}} = \dfrac{\sqrt{3}}{2\sqrt{r^2 - r + 1}}.$$

(ii) From the Law of Sines, $c = \dfrac{b \sin \omega}{\sin \varphi} = \dfrac{b \sin 60}{\sin \varphi} = \dfrac{b}{\sin \varphi} \cdot \dfrac{\sqrt{3}}{2}$

$$\boxed{c = \dfrac{b\sqrt{3(r-1)^2 + (r+1)^2}}{2}} = b\sqrt{r^2 - r + 1}.$$

(iii) When $r = 2$, we obtain from the above formulas $\sin \theta = 1$ and $\sin \varphi = \dfrac{1}{2}$;

$\theta = 90$ and $\varphi = 30$.

(iv) We invite the reader to verify the following claims by using the fact that the quadratic trinomial in $r$, $f(r) = 3(r-1)^2 + (r+1)^2$; $f(r) = 4r^2 - 4r + 4$, has absolute minimum value 3, which occurs at $r = \dfrac{1}{2}$; $f(r) = 4\left(r - \dfrac{1}{2}\right)^2 + 3$.

Thus, on the open interval $(1, +\infty)$, the trinomial $f(r)$ is increasing in value. Accordingly, as $r$ is decreasing in value toward 1, the angle measure $\varphi$ decreases toward the value of 60; and since $\theta = 180 - \omega - \varphi = 180 - 60 - \varphi = 120 - \varphi$, the degree measure $\theta$ is increasingly approaching the value 60, so the limiting position of the triangle at hand, is an equilateral one. On the other hand, as $r$ increases in value toward 2, the angle measure $\varphi$ decreases toward 30, while $\theta$ increasingly approaches 90. Of course, when $r = 2$, we have the exact values $\varphi = 30$ and $\theta = 90$. As $r$ increases away from 2 toward $+\infty$, $\varphi$ decreases toward 0, while $\theta$ increases toward 120.

**References**

[1] Robert Blitzer, *Precalculus Essentials, 2$^{nd}$ Edition*, Prentice-Hall, Inc., 2007, pp. 639-643, 708 pp.